\magnification=\magstep1
\overfullrule=0pt
\def\proofmode{\baselineskip=18pt}
\font\caps=cmcsc10
\font\titlefont=cmssbx10 at 12pt
\font\cbf=cmssbx10
\font\ncbf=cmssbx10 at 9pt
\def\title#1{\topinsert\vskip.5in\endinsert\centerline{\titlefont#1}\bigskip}
\def\author#1{\centerline{\caps#1}\medskip}
\def\abstract#1{{\ninepoint{\narrower\smallskip\noindent
	{\ncbf Abstract.} #1\smallskip}\vskip.5truein}}
	\outer\def\isection#1\par{\vskip0pt plus.3\vsize\vskip0pt plus-.3\vsize
	\bigskip\vskip\parskip\leftline{\cbf#1}\nobreak\smallskip\indent}
\outer\def\proclaim #1. #2\par{\bigbreak\noindent{\cbf#1.\enspace}{\sl#2}\par
	\ifdim\lastskip<\bigskipamount\removelastskip\penalty55\bigskip\fi}
\outer\def\demo #1. #2\par{\medbreak\noindent{\it#1.\enspace}
	{\rm#2}\par\ifdim\lastskip<\medskipamount\removelastskip
	\penalty55\medskip\fi}
\def\proof{\noindent{\it Proof.}\quad}
\def\qed{~\hfill~\blackbox\medskip}
\def\blackbox{\hbox{\vrule width6pt height7pt depth1pt}}
\def\frac#1#2{{\textstyle{#1\over#2}}}
\def\iitem{\itemitem}
\def\myskip{\noalign{\vskip6pt}}
\def\hangbox to #1 #2{\vskip1pt\hangindent #1\noindent \hbox to #1{#2}$\!\!$}

\def\supp{\mathop{\rm supp}\nolimits}
\def\sgn{\mathop{\rm sign}\nolimits}
\def\varep{\varepsilon}
\def\real{\mathop{{\rm I}\kern-.2em{\rm R}}\nolimits}
\def\nat {\mathop{{\rm I}\kern-.2em{\rm N}}\nolimits}
\newskip\ttglue
\def\ninepoint{\def\rm{\fam0\ninerm}
	\textfont0=\ninerm \scriptfont0=\sixrm \scriptscriptfont0=\fiverm
	\textfont1=\ninei  \scriptfont1=\sixi  \scriptscriptfont1=\fivei
	\textfont2=\ninesy  \scriptfont2=\sixsy  \scriptscriptfont2=\fivesy
	\textfont3=\tenex  \scriptfont3=\tenex  \scriptscriptfont3=\tenex
	\textfont\itfam=\nineit  \def\it{\fam\itfam\nineit}
	\textfont\slfam=\ninesl  \def\sl{\fam\slfam\ninesl}
	\textfont\ttfam=\ninett  \def\tt{\fam\ttfam\ninett}
	\textfont\bffam=\ninebf  \scriptfont\bffam=\sixbf
	\scriptscriptfont\bffam=\fivebf  \def\bf{\fam\bffam\ninebf}
	\tt  \ttglue=.5em plus.25em minus.15em
	\normalbaselineskip=11pt
	\setbox\strutbox=\hbox{\vrule height8pt depth3pt width0pt}
	\let\sc=\sevenrm  \let\big=\ninebig \normalbaselines\rm}
	\font\ninerm=cmr9 \font\sixrm=cmr6 \font\fiverm=cmr5
	\font\ninei=cmmi9  \font\sixi=cmmi6   \font\fivei=cmmi5
	\font\ninesy=cmsy9  \font\sixsy=cmsy6 \font\fivesy=cmsy5
	\font\nineit=cmti9  \font\ninesl=cmsl9  \font\ninett=cmtt9
	\font\ninebf=cmbx9  \font\sixbf=cmbx6 \font\fivebf=cmbx5
	\def\ninebig#1{{\hbox{$\textfont0=\tenrm\textfont2=\tensy
	\left#1\vbox to7.25pt{}\right.$}}}
\title{The Distortion Problem}
\author{E. Odell\footnote{}{\rm E.~Odell was partially supported by NSF Grants  
DMS-8903197, DMS-9208482 and TARP 235.
Th. Schlumprecht was partially supported by NSF Grant DMS 9203753 and LEQSF.}  
and Th. Schlumprecht}
\vskip.2in 
\abstract{We prove that Hilbert space is distortable and, in fact, 
arbitrarily distortable. This means that for all $\lambda >1$ there 
exists an equivalent norm $|\cdot|$ on $\ell_2$ such that for all infinite 
dimensional subspaces $Y$ of $\ell_2$ there exist $x,y\in Y$ with 
$\|x\|_2 = \|y\|_2 =1$ yet $|x| >\lambda |y|$. 
\vskip1pt
We also prove that if $X$ is any infinite dimensional Banach space with 
an unconditional basis then the unit sphere of $X$ and the unit sphere 
of $\ell_1$ are uniformly homeomorphic if and only if $X$ does not contain 
$\ell_\infty^n$'s uniformly.} 

\proofmode 
\isection{1. Introduction}

An infinite dimensional Banach space $X$ is {\it distortable\/} 
if there exists an equivalent norm $|\cdot|$ on $X$ and $\lambda >1$ 
such that for all infinite dimensional subspaces $Y$ of $X$, 
$$\sup \left\{ {|y|\over |z|} : y,z\in S(Y;\|\cdot\|)\right\} >\lambda\ , 
\eqno(1.1)$$ 
where $S(Y;\|\cdot\|)$ is the unit sphere of $Y$. 
R.C.~James [11] proved that $\ell_1$ and $c_0$ are not distortable. 
In this paper we prove that $\ell_2$ is distortable. In fact we shall 
prove that 
$\ell_2$ is {\it arbitrarily distortable\/} 
(for every $\lambda >1$ there exists an equivalent norm on $\ell_2$ 
satisfying (1.1)). 

The distortion problem  is related to stability problems for a wider class 
of functions than the class of equivalent norms. A function $f:S(X)\to\real$ 
is {\it oscillation stable\/} on $X$ if for all subspaces $Y$ of $X$ 
and for all $\varep >0$ there exists a subspace $Z$ of $Y$ with 
$$\sup \bigl\{ |f(y) - f(z)| : y,z\in S(Z)\bigr\} <\varep\ . 
\eqno(1.2)$$
(By {\it subspace\/} we shall mean a closed infinite dimensional linear 
subspace unless otherwise specified.) 
It was proved by V.~Milman (see e.g., [28, p.6] or [26,27]   
that every Lipschitz (or even uniformly continuous) 
function $f:S(X)\to\real$ is {\it finitely\/} 
oscillation stable (a subspace $Z$ of arbitrary finite dimension can 
be found satisfying (1.2)). V.~Milman proved in his fundamental paper [26] 
that if all Lipschitz functions on every unit sphere of every 
Banach space were oscillation 
stable, then every $X$ would isomorphically contain $c_0$ or $\ell_p$ for 
some $1\le p<\infty$. Of course Tsirelson's famous example [38]  dashed 
such hopes and caused Milman's paper to be overlooked. However Milman's work 
implicitly contains the result, rediscovered in [10], that if $X$ does 
not contain $c_0$ or $\ell_p$ $(1\le p<\infty)$ then some subspace of $X$ 
admits a distorted norm. Thus the general distortion problem 
(does a given $X$ contain a distortable subspace?) reduces 
to the case $X= \ell_p$ $(1<p<\infty)$. 

For a given space $X$, every Lipschitz function $f:S(X) \to\real$ is 
oscillation stable if and only if every uniformly continuous 
$g:S(X)\to\real$ is oscillation stable. Indeed if such a $g$ were 
not oscillation stable then there exist a subspace $Y$ of $X$ and reals 
$a<b$ such that 
$$C = \bigl\{ y\in S(Y): g(y) <a\bigr\}\quad\hbox{and}\quad 
D= \bigl\{ y\in S(Y):g(y)>b\bigr\}$$ 
are both asymptotic for $Y$ ($C$ is {\it asymptotic\/} for $Y$ if 
$C_\varep\cap S(Z)\ne \emptyset$ for all subspaces $Z$ of $Y$ 
and all $\varep >0$ where $C_\varep =\{ x: d(C,x)<\varep\}$). Since $g$ is 
uniformly continuous, $d(C,D) \equiv \inf \{\|c-d\| :c\in C$, $d\in D\} >0$ 
and so $f(x) \equiv d(C,x)$ is a Lipschitz function on $S(X)$ that 
does not stabilize. 

If $C$ and $D$ are asymptotic sets for a uniformly convex space $X$ 
with $d(C,D)>0$ then $X$ contains a distortable subspace. For example the norm 
$|\cdot|$ on $X$ whose unit ball is the closed convex hull of 
$(A\cup -A\cup \delta\, Ba\, X)$ 
is  a distortion of a subspace 
for sufficiently small $\delta$ and some choice $A\in \{C,D\}$. 
If $X= c_0$ or $\ell_p$ 
$(1\le p<\infty)$, then by the minimality of $X$ one obtains that every  
uniformly continuous $f:S(X)\to\real$ is oscillation stable if and only 
if $S(X)$ does not contain two asymptotic sets a positive distance apart. 
If $X= \ell_p$ $(1<p<\infty)$ then this is, in turn, equivalent to 
$X$ is not distortable. 

T.~Gowers [7] proved that every uniformly continuous function 
$f:S(c_0)\to \real$ is oscillation stable. Every uniformly continuous 
$f:S(\ell_1)\to\real$ is oscillation stable if and only if $\ell_2$ 
(equivalently $\ell_p$, $1<p<\infty$) is not distortable. 
This is seen by considering the Mazur map [25] 
$M:S(\ell_1) \to S(\ell_2)$ given by 
$M(x_i)_{i=1}^\infty = ((\sgn x_i)\sqrt{|x_i|})_{i=1}^\infty$. 
$M$ is a uniform homeomorphism between the two unit spheres (see e.g., [32], 
lemma~1). 
Moreover, since $M$ preserves subspaces spanned by block bases of the 
respective unit vector bases of $\ell_1$ and $\ell_2$, 
$C$ is an asymptotic set for $\ell_1$ if and only if 
$M(C)$ is an asymptotic set for $\ell_2$. 

Gowers theorem combined with our main result and that of Milman's 
yields the 

\proclaim Theorem 1.1. 
Let $X$ be an infinite dimensional Banach space. Then every 
Lipschitz function $f:S(X) \to \real$ is oscillation stable if and 
only if $X$ is $c_0$-saturated. 

($X$ is {\it $c_0$-saturated\/} if every subspace of $X$ contains 
an isomorph of $c_0$.) 

In Section 2 we consider a generalization of the Mazur map. 
The Mazur map satisfies  for $h= (h_i) \in S(\ell_1)^+$ with  
$h$ finitely supported, $M(h)=x$ where $x\in S(\ell_2)^+$ maximizes 
$E(h,y) \equiv \sum_i h_i\log y_i$ over $S(\ell_2)^+$. 
Furthermore in this case $h= x^*\circ x$ where $x^*$ is the unique support 
functional of $x$ and ``$\circ$'' denotes pointwise multiplication 
of the sequences $x$ and $x^*$. 
These facts are well known. We give a proof in Proposition~2.5. 

The generalization is given as follows. Let $X$ have a 
{\it $1$-unconditional\/} 
normalized basis $(e_i)$. 
This just means that $\big\|\,|x|\,\big\| = \|x\|$ for 
all $x=\sum a_ie_i\in X$ where $|x| = \sum |a_i|e_i$. 

We regard $X$ as a discrete lattice. $c_{00}$ denotes the linear space of 
finitely supported sequences on $\nat$. Thus $X\cap c_{00} = \{x\in X: 
\supp x$ is finite$\}$ where ${\supp (\sum a_ie_i) = \{ i:a_i\ne 0\}}$. 
For $B\subseteq \nat$ and  $x =\sum x_ie_i\in X$ we set $Bx = \sum_{i\in B} 
x_ie_i$. We often write $x= (x_i)$. $\ell_1$ is a particular instance of 
such an $X$ and we use the same notational conventions for $\ell_1$. 

The generalization  $F_X$ of the Mazur  map is defined in terms of an 
auxilliary  map, the {\it entropy\/} function $E: (\ell_1\cap c_{00}) \times 
X\to [-\infty,\infty)$ given by $E(h,x) \equiv E(|h|,|x|) \equiv 
\sum_i |h_i| \log |x_i|$ where $h= (h_i)\in \ell_1\cap c_{00}$ and 
$x= (x_i)\in X$ under the convention $0\log0\equiv0$. Fix $h\in \ell_1\cap 
c_{00}$ and $B= \supp h$. Then there exists a unique $x= (x_i)\in S(X)$ 
satisfying 
\smallskip
\iitem{i)} $E(h,x) \ge E(h,y)$ for all $y\in S(X)$ 
\iitem{ii)} $\supp h = \supp x= B$
\iitem{iii)} $\sgn x_i = \sgn h_i$ for $i\in B$.
\smallskip
\noindent This unique $x$ we denote by $F_X (h)$ and we set $E_X(h) = 
E(h,F_X (h)) = \max \{E(h,y):y\in S(X)\}$. 

Indeed the function 
$E(h,\cdot) : \{x\in S(X)^+ : \supp x\subseteq B\}\to [-\infty,0]$ is 
continuous taking real values on those $x$'s with $\supp x= B$ and taking the 
value $-\infty$ otherwise. Thus there exists $x\in S(X)^+$ satisfying ii) 
and $E(h,x) \ge E(h,y)$  if $y\in S(X)^+$, $\supp y\subseteq B$. 
Since $(e_i)$ is 1-unconditional and $E(h,y) = E(h,By)$ for all $y\in X$, 
we obtain i). 
iii) is then achieved by changing the signs of $x_i$ as needed. The 
uniqueness of $x$ follows from the strict concavity of the log function. 
If $\supp x= \supp y =B$ and $x\ne y$ then $E(h,\frac12 (|x|+|y|)) > \frac12 
E(h,|x|) +\frac12 E(h,|y|)$. 

We discovered the map $E$  in a paper of Gillespie  [6] 
and we thank L.~Weis for bringing that paper to our attention. 
A similar map is considered in [37]. As noted there other authors have also 
worked with this map in various contexts ([20,21], [13], [30], [36], [14]). 
The central objective of some of these earlier papers was to show that elements 
of $S(\ell_1)$ could be written as $x^*\circ x$ with $\|x^*\|=\|x\|=1$. 
Our additional 
focal point is the map $F_X$ itself.  For certain $X$, $F_X$  is 
uniformly continuous. In general  
$F_X$ is not uniformly continuous, but retains enough 
structure (Proposition~2.3) to be extremely useful in Section~3. In addition  
it is known (e.g., [37], lemma~39.3) that whenever $x= F_X(h)$ 
there exists $x^*\in S(X^*)$ with $x^*\circ x=h$. 

We prove (Theorem~3.1) that if $X$ has an unconditional basis and if 
$X$ does not contain $\ell_\infty^n$ uniformly in $n$, then there 
exists a uniform homeomorphism $F: S(\ell_1)\to S(X)$. We prove this 
by reducing the problem, this follows easily from the work of [5] and [23], 
to the case where $X$ has a 
1-unconditional basis and is $q$-concave with constant $1$ for some 
$q<\infty$.  
$X$ is {\it $q$-concave\/} with constant $M_q(X)$ if 
$$\biggl( \sum_{i=1}^n \|x^i\|^q\biggr)^{1/q} 
\le M_q (X) \Big\| \biggl( \sum_{i=1}^n |x^i|^q\biggr)^{1/q}\Big\| 
\eqno(1.3)$$ 
whenever $(x^i)_{i=1}^n \subseteq X$. The vector on the right side of (1.3) 
is computed coordinatewise with respect to $(e_j)$. 
In this particular case the uniform homeomorphism $F$ is the map $F_X$ 
described above. 

One  way to attack the distortion problem is to find a 
distortable space $X$ with a 1-unconditional basis and having say $M_2(X)=1$ 
and possessing a describable pair of separated asymptotic sets. Then use 
the map $F_X$ to pull these sets back to a separated pair (easy) of  
asymptotic sets (not easy) in $S(\ell_1)$. Our original proof that 
$\ell_2$ is distortable was a variation of this idea using $X=T_2^*$, 
the dual of convexified Tsirelson space. However much more is possible 
as was shown to us by B.~Maurey. Maurey's elegant argument is given in 
Section~3 (Theorem~3.4). 
We thank him for permitting us to include it in this paper. 

In Section 3 
we use the map $F_X$ for $X=S^*$, the dual space of the arbitrarily 
distortable space constructed in [34] (see also [35]). 
As shown in [9] and implicitly in [34,35] this space contains 
a sequence of nearly biorthogonal sets: $A_k \subseteq S(S)$, 
$A_k^* \subseteq Ba(S^*)$ with $A_k$ asymptotic in $S$ for all $k$. 
By ``nearly biorthogonal'' we mean  that for some sequence 
$\varep_i\downarrow 0$, $|x_k^*(x_j)| <\varep_{\min (k,j)}$ if $k\ne j$, 
$x_k^* \in A_k^*$, $x_j\in A_j$, 
and $A_k^*$ $1-\varep_k$-norms $A_k$. 
The latter means that for all $x_k\in A_k$ there exists $x_k^*\in A_k^*$ 
with $x_k^* (x_k)>1-\varep_k$. 
The particular description of these sets 
is used along with the mapping $F_{S^*}$ to show that the sets 
$$\eqalign{
C_k&\equiv \bigl\{ x\in \ell_2 : |x| = \sqrt{|x_k^*\circ x_k|/\|x_k^*
\circ x_k\|_1 }\ \hbox{ for some }\cr
&\qquad 
\ x_k^*\in A_k^*\ ,\ x_k\in A_k\ \hbox{ with }\ \|x_k^* \circ x_k\|_1\ge 
1-\varep_k \bigr\}\cr}$$
are nearly biorthogonal in $\ell_2$ (easy) and that $C_k$ is asymptotic 
in $\ell_2$.  By $x^* \circ x$ we mean again 
the element of $\ell_1$ given by the 
operation of pointwise multiplication. Thus if $x^* = \sum a_i e_i^*$ 
and $x= \sum b_i e_i$, $x^*\circ x = (a_ib_i)_{i=1}^\infty$. 
$\|\cdot\|_1$ is the $\ell_1$-norm. 

The sets $C_k$ easily lead to an arbitrary distortion of $\ell_2$. In fact 
using an argument of [9] one can prove the following (see also 
Theorem~3.1). 

\proclaim Theorem 1.2. 
For all $1<p<\infty$, $\varep >0$  and $n\in \nat$ 
there exists an equivalent norm $|\cdot|$ on $\ell_p$ such 
that for any block basis $(y_i)$ of the unit vector basis of 
$\ell_p$ there exists a finite block 
basis $(z_i)_{i=1}^n$ of $(y_i)$  
which is $1+\varep$-equivalent to the first $n$ 
terms of the summing basis, $(s_i)_{i=1}^n$. 

The summing basis norm is  
$$\Big\| \sum_{i=1}^n a_i s_i\Big\| 
= \sup \biggl\{ \Big| \sum_{i=1}^\ell a_i\Big| : \ell\le n\biggr\}\ .$$
Thus for all $\lambda >1$ there exists an equivalent norm 
$|\cdot|$ on $\ell_p$ such that no basic sequence in $\ell_p$ is 
$\lambda$-unconditional in the $|\cdot|$ norm. The sets $C_k$, in addition 
to being nearly biorthogonal, are unconditional and spreading  (defined 
in Section~3 just before the statement of Theorem~3.4) and 
seem likely to prove useful elsewhere. 

T.~Gowers  [8] proved the 
conditional theorem that if every equivalent norm on $\ell_2$ admits an almost 
symmetric subspace, then $\ell_2$ is not distortable. Theorem~1.2 shows 
that one cannot even obtain an almost $1$-unconditional subspace in general. 
An easy consequence [10] is that there exists an 
asymptotic set $C\subseteq S(\ell_2)$ 
with ${d(C,-C)>0}$. 

The paper by Lindenstrauss and Pe{\l}czy\'nski [17] also contains some 
nice results on distortion. They consider a restricted form of distortion 
in which the subspace $Y$ of (1.1) is isomorphic to $X$. 

Our notation is standard Banach space terminology as may be found in the 
books [18] and [19]. In Section~2 we use a number of 
results in [5] although 
we cite the corresponding statements in [19]. 

Thanks are due to numerous people, especially B.~Maurey and 
N.~Tomczak-Jaegermann. As we noted, Maurey gave us the elegant argument 
of Section~3. The idea of exploiting the ramifications of being able to 
write elements of $S(\ell_2)$ 
as $\sqrt{x^*\circ x}$ with $x$ in the sphere of a Tsirelson-type 
space $X$ and $x^*\in S(X^*)$ 
in attacking the 
distortion problem is due to Tomczak-Jaegermann. 

\isection{2. Uniform homeomorphisms between unit spheres}

The main result of this section is 

\proclaim Theorem 2.1. 
Let $X$ be a Banach space with an unconditional basis. Then $S(X)$ and 
$S(\ell_1)$ are uniformly homeomorphic if and only if $X$ does not 
contain $\ell_\infty^n$ uniformly in $n$. 

A {\it uniform homeomorphism\/} between two metric spaces is an invertible 
map such that both the map and its inverse are uniformly continuous. 
Many results are known concerning uniform homeomorphisms between 
Banach spaces (see [1] for a nice survey of these results). Our focus 
however is on the unit spheres of Banach spaces. The prototype of such maps  
is the Mazur map discussed in the introduction. 

Before proceeding we set some notation. Unless stated otherwise $X$ 
shall be a Banach space with a normalized $1$-unconditional basis $(e_i)$. 
We regard $X$ as a discrete lattice. $x=(x_i)\in X$ means that 
$x=\sum x_i e_i$, $|x| = (|x_i|)$, 
and $Ba (X)^+ = \{x\in Ba (X): x=|x|\}$. 
$Ba(X)$ is the closed unit ball of $X$. 
For $1\le p<\infty$, $X$ is {\it $p$-convex\/} with $p$-convexity constant 
$M^p(X)$ if for all $(x^i)_{i=1}^n\subseteq X$, 
$$\Big\| \biggl( \sum_{i=1}^n |x^i|^p\biggr)^{1/p}\Big\| 
\le M^p (X) \biggl( \sum_{i=1}^n \|x^i\|^p\biggr)^{1/p}\ ,$$
where $M^p(X)$ is the smallest constant satisfying the inequality. 
The {\it $p$-convexification\/} of $X$ is the Banach space given by 
$$X^{(p)} = \biggl\{ (x_i) : \|(x_i)\|_{(p)} \equiv 
\Big\| \sum_i |x_i|^p e_i\Big\|^{1/p} <\infty \biggr\}\ .$$ 
The unit vector basis of $X^{(p)}$, which we still denote by $(e_i)$, 
is a $1$-unconditional basis for $X^{(p)}$ and $M^p (X^{(p)})=1$. These 
facts may be found in [19, section~1.d].

Let $F_X :\ell_1\cap c_{00}\to S(X)$ be as defined in the introduction. 
As we shall see in Proposition 2.5, $F_X$ generalizes the Mazur map. 
If $X= \ell_p$ ($1<p<\infty$) and $h\in S(\ell_1)^+\cap c_{00}$ then 
$F_X(h) = (h_i^{1/p})$. Even in this nice setting however we cannot use 
our definitions directly on infinitely supported elements. Indeed one can 
find $h\in S(\ell_1)$ with $E_{\ell_2} (h) = -\infty$. The map $F_{\ell_2}$ 
is uniformly continuous on $S(\ell_1)\cap c_{00}$, though, and thus extends 
to a map on $S(\ell_1)$. $E_X$ is not uniformly continuous on 
$S(\ell_1)\cap c_{00}$ but has some positive features as 
the next proposition reveals. 
Some of our arguments could be slightly shortened by referring to the 
papers [20,21], [13], [37] and [6] but we choose to present complete proofs. 

First we define a function 
$\psi(\varep)$ that appears in Proposition~2.3. 
Note that 
there exists a function $\eta :(0,1)\to (0,1)$ so that 
$$\log \frac12 \left(\sqrt{a} +{1\over\sqrt a}\right)>\eta (\varep)
\ \hbox{ if }\ |a-1|>\varep\ \hbox{ with }\ a>0\ .  
\eqno(2.1)$$ 
Indeed, let  $g(a) = \log\frac12 (a+\frac1a)$ for $a>0$. $g$ is continuous on 
$(0,\infty)$, strictly decreasing on $(0,1)$ and strictly increasing on 
$(1,\infty)$. The minimum value of $g$ is $g(1)=0$. Thus there exists 
$\eta :(0,1)\to (0,1)$ so that $|a-1| >\varep$ implies 
$g(\sqrt{a}\,) > \eta (\varep)$.\qed 

\proclaim Definition 2.2. 
$\psi(\varep) = \varep \eta (\varep)$ for $\varep \in (0,1)$. 

\proclaim Proposition 2.3. 
Let $X$ have a $1$-unconditional basis. 
\vskip1pt
\iitem{\rm A.} 
Let $h\in S(\ell_1)^+\cap c_{00}$, let $\varep>0$ and $v\in Ba (X)^+$ be 
such that $E(h,v) \ge E_X(h) -\psi(\varep)$. 
Then if $u= F_X(h)$ there exists $A\subseteq \supp h$ 
satisfying $\|Ah\|>1-\varep$ 
and $(1-\varep)Au \le Av\le (1+\varep) Au$ (the latter inequalities being 
pointwise in the lattice sense). 
\vskip1pt
\iitem{\rm B.} 
Let $h_1,h_2 \in S(\ell_1)^+\cap c_{00}$ with $\|h_1-h_2\|\le 1$. 
Let $x_i= F_X(h_i)$ for $i=1,2$. Then  
$$\Big\| {x_1+x_2\over2} \Big\| \ge 1-\sqrt{\|h_1-h_2\|}\ .$$ 
\smallskip

\proof 
A. Let $u= (u_i)$ and $v= (v_i)$ be as in the statement of (A). 
We may assume that $\supp u = \supp v = B\equiv \supp h$. 
$E(h,v) \ge E_X(h) - \psi(\varep)$ yields 
$$\psi(\varep) \ge \sum_{i\in B} h_i (\log u_i - \log v_i) \ .
\eqno(2.2)$$ 
Since ${u+v\over2} \in Ba (X)^+$ and $u= F_X(h)$ from (2.2) we obtain 
$$\eqalign{ 
\psi (\varep) & \ge \sum_{i\in B}h_i\left[ \log \Bigl({u_i+v_i\over2}\Bigr) - 
\log v_i \right] \cr 
\myskip 
&= \sum_{i\in B} h_i \left[ \frac12 \log u_i +\frac12\log v_i +\log 
\Bigl( {u_i+v_i\over2}\Bigr) - \log\sqrt{u_iv_i} - \log v_i\right]\cr 
\myskip 
&= \frac12 \sum_{i\in B} h_i (\log u_i-\log v_i) 
+ \sum_{i\in B} h_i \log \frac12 \Bigl(\sqrt{v_i\over u_i} + 
\sqrt{u_i\over v_i}\ \Bigr) \ .\cr}$$ 
The first term in the last expression is nonnegative so 
$$\psi (\varep) \ge \sum_{i\in B} h_i \log \frac12 \Bigl( \sqrt{v_i\over u_i} 
+ \sqrt{u_i\over v_i}\ \Bigr)\ .
\eqno(2.3)$$ 
Now $|{v_i\over u_i} -1|<\varep$ if and only if $(1-\varep) u_i\le v_i 
\le (1+\varep)u_i$. Let $I= \{i\in B: |{v_i\over u_i}-1|\ge \varep\}$. 
For $i\in I$, 
$$\log \frac12 \Bigl( \sqrt{u_i\over v_i} + \sqrt{v_i\over u_i}\ \Bigr) 
\ge \eta (\varep)\qquad \hbox{(by (2.1))}\ . 
\eqno(2.4)$$ 
Let $J= \{i\in B :\log \frac12 (\sqrt{u_i\over v_i}+\sqrt{v_i\over u_i}\,) 
\ge \eta (\varep)\}$. 
Thus $I\subseteq J$  by (2.4) and from (2.3),  
$$\sum_{i\in I} h_i \le \sum_{i\in J}h_i 
\le {1\over\eta(\varep)} \sum_{i\in J} h_i \log \frac12 
\Bigl(\sqrt{u_i\over v_i} + \sqrt{v_i\over u_i}\ \Bigr) 
\le {\psi (\varep)\over \eta(\varep)} = \varep\ .$$ 
Thus (A) follows with $A= B\setminus I$. 

B. 
Let $\| {x_1+x_2\over2}\| \equiv 1-2\varep$. 
Set $\tilde x_1 = x_1 + \varep x_2$ and 
$\tilde x_2 = x_2 +\varep x_1$. Thus $\supp \tilde x_1 = 
\supp \tilde x_2 = \supp h_1 \cup \supp h_2$ and 
$\| {\tilde x_1+\tilde x_2\over 2}\| \le 1-\varep$. 
We may assume $\varep >0$. For $j\in \supp \tilde x_1$, 
$|\log \tilde x_{1,j} - \log \tilde x_{2,j}| \le |\log \varep|$ 
where $\tilde x_i = (\tilde x_{i,j})$ for $i=1,2$. 

>From this and $\tilde x_1 \ge x_1$ we obtain 
$$\eqalign{E(h_1,\tilde x_1)
&\ge E(h_1,x_1)\ge E\Bigl(h_1,{\tilde x_1+\tilde x_2\over 2(1-\varep)}\Bigr)\cr 
\myskip 
&= E\bigl( h_1, {\tilde x_1+\tilde x_2\over2}\Bigr) +|\log (1-\varep)|\cr 
\myskip 
&\ge \frac12 E(h_1,\tilde x_1) +\frac12 E(h_1,\tilde x_2) 
+ |\log (1-\varep)|\ .\cr}$$ 
Thus 
$$ |\log (1-\varep)| \le \frac12 \bigl( E(h_1,\tilde x_1) - 
E(h_1,\tilde x_2)\bigr)\ .$$
Similarly, 
$$|\log (1-\varep)| \le \frac12 \bigl( E(h_2,\tilde x_2)- 
E(h_2,\tilde x_1)\bigr)\ .$$ 
Averaging the two inequalities yields 
$$\eqalign{\varep &\le |\log (1-\varep)| 
\le\frac14 \bigl( E(h_1,\tilde x_1) - E(h_1,\tilde x_2) -
E(h_2,\tilde x_1) +E(h_2,\tilde x_2) \cr 
\myskip 
&= \frac14 \sum_{j\in B} (h_{1,j} - h_{2,j}) (\log \tilde x_{1,j} 
- \log \tilde x_{2,j}) \cr 
\myskip 
&\le \frac14 \|h_1 - h_2\| \, |\log\varep| 
\le\frac14 \|h_1-h_2\| \varep^{-1}\ .\cr}$$ 
Thus $\varep \le \frac12 \|h_1-h_2\|^{1/2}$. 
Hence $\|{x_1+x_2\over2}\| = 1-2\varep \ge 1-\|h_1-h_2\|^{1/2}$.\qed 

\proclaim Proposition 2.4. 
Let $X$ be a uniformly convex Banach space with a $1$-unconditional 
basis. The map $F_X :S(\ell_1)\cap c_{00}\to S(X)$ is 
uniformly continuous. Moreover the modulus of continuity of $F_X$ 
depends solely on the modulus of uniform convexity of~$X$. 

\proof 
The uniform continuity of $F_X$ on $S(\ell_1)^+\cap c_{00}$ follows 
immediately from Proposition~2.3(B).

Precisely, there is a function $g(\varep)$, depending solely upon the 
modulus of uniform convexity of $X$, which is continuous at $0$ with 
$g(0)=0$ and  satisfies 
$$\| F_X(h_1) - F_X(h_2)\| \le g(\|h_1-h_2\|)$$ 
for $h_1,h_2\in S(\ell_1)^+ \cap c_{00}$. 
A consequence of this is that if $h\in S(\ell_1)^+ \cap c_{00}$, 
$x=F_X(h)$ and $I\subseteq \nat$ is such that $\|Ih\|<\varep$ then 
$\|Ix\|<g(2\varep)$. Indeed if $J=\nat \setminus I$, 
$$\Big\| h- {Jh\over \|Jh\|}\Big\| = \|Ih\| + 
\Big\| Jh - {Jh\over \|Jh\|}\Big\| <2\varep\ .$$ 
Thus since $Ix = I(F_X(h) - F_X (Jh/\|Jh\|))$, 
$$\|Ix\| \le \Big\| F_X(h) - F_X\Bigl( {Jh\over \|Jh\|}\Bigr)\Big\| 
\le g(2\varep)\ .$$ 

For the general case let $h_1,h_2 \in S(\ell_1) \cap c_{00}$ with 
$\|h_1-h_2\| = \varep$. Let $F_X(|h_i|) = |x_i|$ for $i=1,2$. 
Then $x_i\equiv \sgn h_i \circ |x_i|$, ``$\circ$'' denoting 
pointwise multiplication, satisfies $x_i= F_X(h_i)$ for $i=1,2$. 
Also $\|\, |h_1| -|h_2|\,\| \le \|h_1-h_2\|$.  
Thus 
$$\|x_1-x_2\| \le \big\|\,|x_1| - |x_2|\,\big\| 
+ \Big\| \sum_{j\in I} (|x_{1,j}| + |x_{2,j}|) e_j\Big\|$$ 
where $I= \{j: \sgn x_{1,j} \ne \sgn x_{2,j}\}$
$$\eqalignno{
&\le g\bigl(\big\|\,|h_1| - |h_2|\,\big\|\bigr) 
+ \big\| I|x_1|\,\big\| + \big\| I|x_2|\,\big\|\cr 
\myskip 
&\le g(\varep) + g(2\varep) + g(2\varep)\ .&\blackbox\cr}$$ 
\medskip 

Here is a fact we promised earlier. 

\proclaim Proposition 2.5. 
Let $X=\ell_p$, $1<p<\infty$. Then $F_X$ is the Mazur map, i.e., 
if $h\in S(\ell_1)^+ \cap c_{00}$ then $F_X(h) = (h_i^{1/p})$.  

\proof 
Let $h\in S(\ell_1)^+ \cap c_{00}$, $B= \supp h$ and   $F_X(h)=x$. 
Then $\supp x=B$ and the vector $(x_i)_{i\in B}$ maximizes the function 
$\real_+^B \ni (y_i)\mapsto \sum_{i\in B} h_i\log y_i$ under the 
restriction $\sum_{i\in B} y_i^p = 1$. By the method of Lagrange 
multipliers this implies that there is a number $c\ne 0$ so that 
${h_i\over x_i} = cpx_i^{p-1}$ for $i\in B$. 
Thus $x_i = (cp)^{-1/p} h_i^{1/p}$. Since $\|x\|_p=1$, 
$$c= p^{-1}\ \hbox{ and }\ x_i = h_i^{1/p}\ \hbox{ for }\ i\in B\ .
\eqno\blackbox$$
\medskip

If $X$ is uniformly convex, by Proposition 2.4 
the map $F_X$ extends uniquely to a uniformly 
continuous map, which we still denote by $F_X$, from $S(\ell_1)\to S(X)$. 

\proclaim Proposition 2.6.
Let $X$ be a uniformly convex uniformly smooth Banach space with a 
$1$-unconditional basis. Then $F_X : S(\ell_1)\to S(X)$ is invertible 
and $(F_X)^{-1}$ is uniformly continuous, with modulus of continuity 
depending only on the modulus of uniform smoothness of $X$. 
For $x\in S(X)$, $F_X^{-1}(x) = \sgn (x)\circ x^*\circ x=|x^*|\circ x$ 
where $x^*$ is the unique support 
functional of $x$.

\proof 
For $x\in S(X)$ there exists a unique element $x^* \in S(X^*)$ such that 
$x^* (x)=1$. The biorthogonal functionals $(e_i^*)$ are a 
$1$-unconditional basis for $X^*$ and thus we can express $x^* = \sum 
x_i^* e_i^*$ and write $x^* = (x_i^*)$. The element $x^*\circ x 
\in S(\ell_1)^+$ and $\sgn x^*=\sgn x$. Let $G(x) = |x^*| \circ x$. 
$G$ is uniformly continuous. Indeed the map 
$S(X) \ni x\mapsto x^*$, the supporting functional, is uniformly 
continuous since $X$ is uniformly smooth. The modulus of continuity of this 
map depends solely on the modulus of uniform smoothness of $X$ 
(see e.g., [3], p.36). Let $G(x_i) = h_i = |x_i^*|\circ x_i$ for 
$i=1,2$. Then 
$$\eqalign{ \|h_1-h_2\| 
&= \|\,|x_1^*|\circ x_1 - |x_2^*|\circ x_2\|\cr  
&\le \|\,|x_1^*|\circ (x_1 -x_2)\| \cr
&\qquad + \|(|x_1^*| - |x_2^*|) \circ x_2\|\cr
&\le \|x_1^*\| \, \|x_1-x_2\| 
+ \|\,|x_1^*|-|x_2^*|\|\, \|x_2\|\cr 
&\le \|x_1-x_2\|  + \|x_1^* - x_2^*\|\cr}$$
which proves that $G$ is uniformly continuous. 

It remains only to show that $G= F_X^{-1}$. 
Since $G(x) =\sgn x\circ G(|x|)$ we need only show that 
$G(F(h)) = h$ for $h\in S(\ell_1)^+ \cap c_{00}$ and $F(G(x))=x$ 
for $x\in S(X)^+\cap c_{00}$. 

If $h\in S(\ell_1^+)\cap c_{00}$ and $x= F_X(h)$ then, as in the proof of 
Proposition~2.5, the method of Lagrange multipliers yields that 
$\vec\nabla E(h,x) = (h_i/x_i)_{i\in \supp h}$ 
equals a multiple of $(x_i^*)_{i\in \supp h}$ where $x^*$ is the 
support functional of $x$. This multiple must be $1$ and $h_i = x_i^*\circ 
x_i$ or $G(F(h)) =h$. 

That $F(G(x))=x$ follows once we observe that if $h= x^*\circ x=y^*\circ y$, 
all norm 1 elements, then $x=y$. Assume for simplicity $\supp h= \{1,2,
\ldots,n\}$. Define $f(z) = \|z\| - E(h,z)$ for $z\in U$, a convex open 
subset of the positive cone $Ba (\langle e_i\rangle_{i=1}^n)^+$ 
which contains both $x$ and $y$ and is bounded away from the boundary of the 
cone. $f(z)$ is strictly convex so $\vec\nabla f(z) = \vec0$ for at 
most one point. But $\vec\nabla f(z)=\vec0$ iff $h= z^*\circ z$.\qed  

\proclaim Corollary 2.7. {\rm [37, lemma 39.3].} 
Let $X$ have a $1$-unconditional basis and let $h\in S(\ell_1^+)\cap c_{00}$ 
with $x\in F_X(h)$. Then  there exists $x^*\in S(X^*)$ with $x^*\circ x=h$. 

\proof 
We may restrict our attention to $X=\langle e_i\rangle_{i\in \supp h}$. 
The result follows if $X$ is smooth from the proof of Proposition~2.6.  
Let $\|\cdot\|_n$ be a sequence of smooth norms on $X$ with 
$\|\cdot\|_n \to \|\cdot\|$ and such  that ${x\over \|x\|_n} \in F_{X_n}(h)$. 
Then use a compactness argument.\qed 

Before proving  Theorem 2.1 we need one more proposition. 
Recall that $X^{(p)}$ is the $p$-convexification of $X$. The map $G_p$ 
below is another generalization of the Mazur map. 

\proclaim Proposition 2.8. 
Let $1<p<\infty$ and let $X$ be a Banach space with a $1$-unconditional 
basis. The map $G_p: S(X^{(p)})\to S(X)$ given by $G_p (x) = 
\sgn (x)\circ |x|^p = ((\sgn x_i)
|x_i|^p)$ for $x=(x_i)$ is a uniform homeomorphism. 
Moreover the modulus of continuity of $G_p$ and $G_p^{-1}$ are functions 
solely of $p$. 

\proof 
As usual $(e_i)$ denotes the normalized $1$-unconditional basis of both 
$X$ and $X^{(p)}$. 
Let $x,y\in S(X^{(p)})$ with $\delta \equiv \|x-y\|_{(p)}$. We shall show that 
$$2^{1-p} \delta^p \le \|G_p(x) - G_p(y)\| 
\le \delta^p +\delta^{p/2} + 2\bigl( 1-(1-\sqrt{\delta}\,)^p\bigr)$$ 
which will complete the proof. 

Let $x= \sum x_i e_i$ and $y= \sum y_i e_i$. 
$$\eqalign{ \|G_p(x) - G_p(y)\| 
& = \Big\| \sum_{i=1}^\infty \bigl( \sgn (x_i) |x_i|^p - \sgn (y_i)|y_i|^p
\bigr) e_i\Big\| \cr 
\myskip 
& = \Big\| \sum_{i\in I_+} (|x_i|^p - |y_i|^p) e_i 
+ \sum_{i\in I_-} (|x_i|^p + |y_i|^p) e_i\Bigr\|\cr}$$ 
where 
$$I_+ = \bigl\{ i: \sgn (x_i) = \sgn (y_i)\bigr\}\quad\hbox{and}\quad 
I_- = \bigl\{ i:\sgn (x_i) \ne \sgn (y_i)\bigr\}\ .$$ 
We denote the two terms in the last norm expression as $d_+$ and $d_-$, 
respectively.

Since $a^p-b^p \ge (a-b)^p$ and $a^p + b^p \ge 2^{1-p} (a+b)^p$ 
for $a\ge b\ge0$ we deduce from the $1$-unconditionality of $(e_i)$ that 
$$\eqalign{ \|d_+ + d_-\| 
&\ge \Big\| \sum_{i\in I_+} \big|\, |x_i|-|y_i|\,\big|^p e_i 
+ 2^{1-p} \sum_{i\in I_+} (|x_i| + |y_i|)^p e_i\Big\| \cr 
\myskip 
&\ge 2^{1-p}\Big\|\sum |x_i-y_i|^p e_i\Big\|= 2^{1-p} \| x-y\|_{(p)}^p\ .\cr}$$ 

To prove the upper estimate we begin by noting that 
$$\eqalign{ 
\|d_-\| & \le \Big\| \sum_{i\in I_-} |x_i-y_i|^p e_i \Big\| \cr 
\myskip 
&\le \| x-y\|_{(p)}^p = \delta^p\ .\cr}$$ 
Set $q= 1-\sqrt{\delta}$ and $c= (1-q)^{-p} = \delta^{-p/2}$. 
For $a,b\ge 0$ with $0\le b\le qa$ we have 
$$c(a-b)^p - (a^p-b^p) \ge c(1-q)^p a^p-a^p = a^p \bigl(c(1-q)^p-1\bigr) 
=0\ .\eqno(2.5)$$ 
Let $I'_+ = \{i\in I_+ :|y_i| < q|x_i|$ or $|x_i| <q|y_i|$ and 
$I''_+ = I_+\setminus I'_+$. Write $d_+ = d'_+ + d''_+$ where 
$d'_+ = \sum_{i\in I'_+} (|x_i|^p - |y_i|^p) e_i$ and $d''_+  = d_+-d'_+$. 
Thus (2.5) yields that 
$$\eqalign{
\|d'_+\| & \le c\Big\| \sum_{i\in I'_+} \big|\, |x_i| - |y_i|\,\big|^p e_i 
\Big\| \cr 
\myskip 
& \le \delta^{-p/2} \|x-y\|_{(p)}^p = \delta^{p/2}\ .\cr}$$ 
Furthermore, 
$$\eqalignno{ \|d''_+\| 
& \le (1-q^p) \Big\| \sum_{i\in I''_+} (|x_i|^p + |y_i|^p)e_i\Big\| \cr 
\myskip 
& \le 2(1-q^p) \le 2\bigl( 1-(1-\sqrt{\delta}\,)^p\bigr)\ .
&\blackbox\cr}$$

\demo Proof of Theorem 2.1. 
It follows quickly from work of Enflo that if $X$ contains $\ell_\infty^n$ 
uniformly in $n$ then $S(X)$ is not uniformly homeomorphic to a subset of 
$S(\ell_1)$. Indeed Enflo [4] proved that a certain family of finite 
subsets of $Ba(\ell_\infty^n)$, $n\in \nat$, cannot be uniformly embedded into 
$Ba(\ell_2)$ and hence neither into $Ba(\ell_1)$. But 
$B(\ell_\infty^n)$ embeds 
isometrically into $S(\ell_\infty^{n+1})$
and hence these finite subsets embed uniformly into $S(X)$. 
\vskip1pt
For the converse assume that $X$ does not contain $\ell_\infty^n$ 
uniformly in $n$. We may suppose that $X$ has a $1$-unconditional 
basis $(e_i)$. 
\vskip1pt
By a theorem of Maurey and Pisier [23], $X$ has cotype $q'$ for some 
$q'<\infty$. This implies that $X$ is $q$-concave for all $q>q'$ 
([19, p.88]). Fix $q>q'$. There exists an equivalent norm on $X$ for 
which $(e_i)$ is still $1$-unconditional and $M_q (X)=1$ 
([19, p.54]). The 2-convexification of $X$ in this norm, $X^{(2)}$, satisfies 
$M_{2q} (X^{(2)}) = 1= M^2 (X^{(2)})$ ([19, p.54]). In particular 
$X^{(2)}$ is uniformly convex and uniformly smooth ([19, p.80]) 
and so $F_{X^{(2)}} : S(\ell_1) \to S(X^{(2)})$ is a uniform 
homeomorphism by Proposition~2.6. Thus $G_2\circ F_{X^{(2)}} :S(\ell_1) 
\to S(X)$ is a uniform homeomorphism by Proposition~2.8.\qed 

\demo Remark. 
If $X$ has a 1-unconditional basis and $M_q(X) =1$ for some $q<\infty$, 
the map $G_2 \circ F_{X^{(2)}} =F_X$. Furthermore the modulus of 
continuity of $F_X$ and $F_X^{-1}$ are functions solely of $q$. 

The uniform homeomorphism theorem extends to unit balls by the 
following simple proposition. 

\proclaim Proposition 2.9. 
Let $X$ and $Y$ be Banach spaces and let $F:S(X) \to S(Y)$ be a uniform 
homeomorphism. For $x\in Ba (X)$ let $\bar F (x) = \|x\| F(x/\|x\|)$ 
if $x\ne 0$ and $\bar F (0)=0$. Then $\widetilde F$ is a uniform homeomorphism 
between $Ba (X)$ and $Ba (Y)$. 

\proof 
Clearly $\bar F$ is a bijection. Since $\bar F^{\,-1} (y) = \|y\| F^{-1} 
(y/\|y\|)$ for $y\ne 0$, it suffices to show that $\bar F$ is uniformly 
continuous. Let $f$ be the modulus of continuity of $F$, 
i.e., $\|F(x_1) - F(x_2)\| \le f(\|x_1-x_2\|)$. 

Let $x_1,x_2 \in Ba (X)$ with $\|x_1-x_2\| = \delta$, 
$\lambda_1 = \|x_1\|$, $\lambda_2 = \|x_2\|$ and $\lambda_1 \ge \lambda_2$. 
$$\eqalign{ \|\bar F(x_1) - \bar F(x_2)\| 
& = \Big\| \lambda_1 F\Bigl( {x_1\over\lambda_1}\Bigr) 
- \lambda_2 F\Bigl( {x_2\over\lambda_2}\Bigr)\Big\| \cr 
& \le  (\lambda_1-\lambda_2) + \lambda_2 
\Big\| F\Bigl( {x_1\over\lambda_1}\Bigr) 
- F\Bigl( {x_2\over\lambda_2}\Bigr)\Big\|\cr}$$ 
If $\lambda_2 < \delta^{1/4}$ this is less than $\delta + 2\delta^{1/4}$. 
Otherwise 
$$\eqalign{ 
\Big\|{x_1\over\lambda_1} - {x_2\over\lambda_2}\Big\| 
& = {1\over\lambda_1\lambda_2} \|\lambda_2x_1-\lambda_1 x_2\| \cr 
\myskip 
&\le {1\over \lambda_1\lambda_2} \bigl[ \lambda_1 \|x_1-x_2\| 
+ \lambda_1-\lambda_2\bigr] 
\le  {\delta\over \lambda_2} + {\delta\over\lambda_1\lambda_2}\cr 
\myskip 
&\le {2\delta\over \lambda_1\lambda_2} 
\le  {2\delta\over\sqrt\delta} = 2\sqrt\delta\ .\cr} $$ 
Thus 
$$\|\bar F (x_1)-\bar F(x_2)\| \le \max \bigl(\delta + f( 2\sqrt\delta\,)\ ,\ 
\delta +2\delta^{1/4}\bigr)\ .\eqno\blackbox$$ 
\medskip

\demo Remark. 
It is not possible, in general, to replace ``uniformly homeomorphic''
by ``Lipschitz equivalent'' in Theorem~2.1. Indeed if $S(X)$ and $S(Y)$ 
are Lipschitz equivalent, then an argument much like that of 
Proposition~2.9, yields that $X$ and $Y$ are Lipschitz equivalent 
which need not be true (see [1]). 

There exist separable infinite dimensional Banach spaces $X$ not 
containing $\ell_\infty^n$'s uniformly such that $Ba (X)$ does not 
embed uniformly into $\ell_2$. 
For example the James' nonoctohedral space [12] has this property. 
Indeed, Y.~Raynaud [31] proved that 
if $X$ is not reflexive and $Ba(X)$ embeds uniformly into $\ell_2$,  
then $X$ admits an $\ell_1$-spreading model. 

Fouad Chaatit [2] has extended Theorem 2.1. He showed one can replace 
the hypothesis that $X$ has an unconditional basis with the more 
general assumption that $X$ is a separable infinite dimensional 
Banach lattice. 
N.J.~Kalton [15] has subsequently discovered another proof of this result 
using complex interpolation theory.

\isection{3. $\ell_2$ is arbitrarily distortable} 

Let $X$ be a Banach space with a basis $(e_i)$. 
A {\it block subspace\/} of $X$ is any subspace spanned by a block basis 
of $(e_i)$. $X$ is {\it sequentially arbitrarily distortable\/} if there 
exist a sequence of equivalent norms $\|\cdot\|_i$ on $X$ and 
$\varep_i\downarrow 0$ such that: 
$$\eqalign{
&\|\cdot\|_i\le\|\cdot\|\ \hbox{ for all $i$ and for all 
subspaces  $Y$  of  $X$}\cr
&\hbox{and for all $i_0\in \nat$ there exists $y\in S(Y,\|\cdot\|_{i_0})$}\cr
&\hbox{with 
$\|y\|_i\le \varep_{\min (i,i_0)}$ for $i\ne i_0$.}\cr}$$  
In the terminology of [9] this is equivalent to saying that $X$ 
contains an asymptotic biorthogonal system with vanishing constants. 

If $X$ is sequentially arbitrarily distortable then $X$ is arbitrarily 
distortable. Even more can be said however. 

\proclaim Theorem 3.1. 
Let $X$ be a sequentially arbitrarily distortable Banach space with a 
basis $(e_i)$. For all $n\in\nat$ and $\varep >0$ 
there exists an equivalent norm $|\cdot|$ 
on $X$ with the following property.  Let $(y_i)_{i=1}^n$ be a monotone 
basis for an $n$-dimensional Banach space. Then every block basis of 
$(e_i)$ admits a further finite block basis $(x_i)_{i=1}^n$ which is 
$(1+\varep)$-equivalent to $(y_i)_{i=1}^n$. 

The space $S$ of [34] was shown in [9] to be sequentially arbitrarily 
distortable. The argument used to prove Theorem~3.1 is a slight variation 
of an argument which appears in [9]  which, in turn, 
has its origins in [24].  

\demo Proof of Theorem 3.1. 
\vskip1pt
Choose for $n\in\nat$ and $\varep>0$, 
$(B_i)_{i=1}^{k(n)}$ a finite sequence of 
$n$-dimensional Banach spaces, each having a monotone basis, such that 
every monotone basis of length $n$ is $(1+\varep)$-equivalent 
to the basis of some 
$B_i^n$. Let $(w_i)_{i=1}^\infty$ be a normalized monotone basis 
for $W\equiv (\sum_{n,i}B_i^n)_{\ell_2}$ such that the monotone basis 
of each $B_i^n$ is 1-equivalent to $(w_i)_{i\in A_i^n}$ for some segment 
$A_i^n\subseteq \nat$. Let $(w_i^*)$ be the biorthogonal functionals 
of $(w_i)$.
\vskip1pt
It suffices to prove that for all $n\in\nat$ there exists an equivalent norm 
$|\cdot|$ on $X$ such that every block basis of $(e_i)$ admits a further 
block basis $(x_i)_{i=1}^n$ which is $(1+{8\over n})$-equivalent 
to $(w_i)_{i=1}^n$.
\vskip1pt 
Let $n\in\nat$, $\varep_i\downarrow 0$ and let $\|\cdot\|_i$ be a 
sequence of equivalent norms on $X$ satisfying the definition of 
sequentially arbitrarily distortable.  Let $\varep >0$ with $n^5\varep <1$. 
We may assume that $\max_i \varep_i < \varep/4$. 
\vskip1pt
Let $X_i = (X,\|\cdot\|_i)$. 
Let $(z_i^*)_{i=2}^\infty$ be an enumeration of all elements of the linear 
span of $(e_i^*)$ which have rational coordinates. Set 
$$\eqalign{
\Gamma =& \biggl\{ z^* = \sum_{i=1}^n b_i  \sum_{j=(i-1)n+1}^{in} z_{k_j}^* 
: k_1<\cdots < k_{n^2}\ ,\ (z_{k_i}^*)_{i=1}^{n^2}\ \hbox{ is a finite }\cr
&\qquad \hbox{block basis of $(e_i^*)$ with }\ z_{k_1}^* \in 3Ba(X_1^*)\ ,\cr
&\qquad z_{k_{i+1}}^* \in 3Ba(X_{k_i}^*)\ \hbox{ for }\ 
1\le i\le n^2-1\hbox{ and } \sum_{i=1}^n b_iw_i^* \in Ba(W^*)\biggr\}\ .\cr}$$
Define $|\cdot|$ on $X$ by 
$$|x| =\sup \bigl\{ | z^* (x) | : z^* \in \Gamma\bigr\}\ .$$ 
Then  $3\|x\|_1 \le |x| \le 6n^2\|x\|$ for all $x\in X$ and so 
$|\cdot|$ is an equivalent norm on $X$. 
\vskip1pt 
Let $Z$ be any block subspace of $X$. Since $X$ cannot contain $\ell_1$,  
we may assume by [33] that $Z$ is spanned by a normalized weakly null block 
basis of $(e_i)$, denoted $(z_i)$. Using the argument that a subsequence 
of $(z_i)$ is nearly monotone for any given norm $|\cdot|_i$ and a 
diagonal argument we may suppose that for all $i$, $\|P_A\|_i < 2.5$ 
whenever $A\subseteq \nat$ is a {\it segment\/} of $\nat$ with 
$i\le \min A$. (Here $P_A$ is the projection $P_A(\sum a_i z_i) 
= \sum_{i\in A} a_i z_i$.) 
\vskip1pt
>From our hypotheses we can then choose block bases $(\bar x_i)_{i=1}^{n^2}$ of 
$(z_i)$, and $(z_{k_i}^*)_{i=1}^{n^2} $ of $(e_i^*)$ satisfying 
$k_1<k_2<\cdots < k_{n^2}$ and 
\vskip1pt
\iitem{i)} $z_{k_1}^* \in 3Ba (X_1^*)$ and $z_{k_{i+1}}^* \in 3Ba (X_{k_i}^*)$ 
for $1\le i<n^2$. 
\iitem{ii)} $z_{k_i}^* (\bar x_j) =\delta_{ij}$ for $1\le i,j\le n^2$. 
\iitem{iii)} $\|\bar x_i\|_j < \varep/3$ if $j\ne k_i$. 
\smallskip
Let 
$x_i = {1\over n} \sum_{j=(i-1)n+1}^{in} \bar x_j$  for $1\le i\le  n$.
\vskip1pt
Let $\|\sum_1^n a_iw_i\|=1= \sum_1^n a_ib_i$ where 
$\|\sum_1^n b_iw_i^*\|=1$. Let $z^* = \sum_{i=1}^n b_i 
\sum_{j=(i-1)n+1}^{in} z_{k_j}^*$ and 
note that $z^*\in\Gamma$. Thus 
$$\Big| \sum_1^n a_i x_i\Big| \ge z^* \biggl( \sum_1^n a_ix_i\biggr) 
= \sum_1^n a_ib_i = 1\ .$$ 
\vskip1pt
For the reverse inequality, let 
$\bar z^* = \sum_{i=1}^n c_i \sum_{j=(i-1)n+1}^{in} z_{m_j}^* \in\Gamma$ with 
$z_{m_1}^* \in 3Ba (X_1^*)$, $z_{m_{i+1}}^* \in 3Ba (X_{m_i}^*)$ for $i<n^2$ 
and $\|\sum_1^n c_i w_i^*\| \le 1$. 
Let $j_0$ be the 
smallest integer such that $m_{j_0} \ne k_{j_0}$. 
We first deduce from the definition of $\Gamma$ and the choice of $(\bar x_i)$ 
that $|z_{m_i}^* (\bar x_j)|<\varep$ and $|z_{m_j}^*(\bar x_i)|<\varep$ 
if $i<j_0$, $j\le n^2$ 
and $i\ne j$. Secondly we claim that 
$$\{ m_{j_0},m_{j_0+1},\ldots,m_{n^2}\} \cap 
\{k_{j_0} ,k_{j_{0+1}},\ldots,k_{n^2}\} = \emptyset\ .$$ 
Indeed if not let $j\ge j_0$ be the smallest integer such that $m_j = k_i$ 
for some $i\ge j_0$. If $j= j_0$ then $i>j_0$. But then (letting $k_0\equiv1$) 
$z_{m_j}^* \in 3Ba (X_{k_{j_0-1}}^*)$ and $\|\bar x_i\|_{k_{j_0-1}} <\varep/3$ 
which contradicts $z_{k_i}^* (\bar x_i)=1$. 
If $j>j_0$ then $z_{m_j}^* \in 3Ba(X_{m_{j-1}}^*)$ and $\|\bar x_i\|_{m_{j-1}}
<\varep/3$ since $m_{j-1} \ne k_{i-1}$, yielding again  a contradiction 
to $z_{k_i}^* (\bar x_i) =1$. 
\vskip1pt 
It follows that $|z_{m_{j_0}}^* (\bar x_i)| <\varep$ if $i\ne j_0$ and 
$|z_{m_j}^* (\bar x_i)| <\varep$ if $j>j_0$ and $i\le n^2$. 
Let $j_0 = i_0n +s_0$ with $0\le i_0<n$, $1\le s_0\le n$. 
Then 
$$\eqalign{
\Big| \bar z^* \biggl( \sum_{i=1}^n a_i x_i\biggr)\Big|  
&= \Big| \biggl( \sum_{i=1}^n c_i \sum_{j=(i-1)n+1}^{in} z_{m_j}^*\biggr) 
\biggl( \sum_{i=1}^n a_i\; {1\over n} \sum_{j=(i-1)n+1}^{in} \bar x_j
\biggr) \Big|\cr 
\myskip 
&\le \Big| \sum_{i=1}^{i_0} c_i a_i + {s_0-1\over n} c_{i_0+1} a_{i_0+1}\Big| 
+ 3\left| {c_{i_0+1} a_{i_0+1} \over n}\right| + n^4 \varep 
\max_i |a_i c_i|\cr 
\myskip 
&\le\Big\|\sum_1^n a_i w_i\Big\| \left[1+{6\over n}+{2\over n}\right]\ .\cr} $$
We used that from monotonicity the first term in the next to last 
inequality does not exceed 
$$\max \biggl(\big|\sum_{i=1}^{i_0} c_i a_i\big| , 
\big|\sum_{i=1}^{i_0+1} c_ia_i\big|\biggr)
\le \Big\|\sum_1^n a_i w_i\Big\|$$ 
and $|c_ia_i| \le 2$ for all $i$.\qed  

\demo Remark. 
The proof of Theorem 3.1 requires only the following condition. 
For all $\varep >0$ there exists a sequence of equivalent norms 
$\|\cdot\|_i\le \|\cdot\|$ on $X$ such that for all subspaces 
$Z$ of $X$ and all $i_0\in \nat$ there exists $y\in S(Z,\|\cdot\|_{i_0})$ 
with $\|y\|_i <\varep$ if $i\ne i_0$.  
In the terminology of [9] this says that for all $\varep>0$, $X$ contains 
an asymptotic biorthogonal system with constant $\varep$. 
Theorem~1.2 is a special case of Theorem 3.1.

Theorem 1.2 yields that a sequentially arbitrarily distortable 
Banach space can be renormed to not contain an almost bimonotone 
basic sequence. Since $\|s_1-2s_2\|=1$, the best constant that can 
be achieved for the norm of the tail projections of a basic sequence is~2. 

Other curious norms can be put on sequentially arbitrarily distortable 
spaces $X$. 
For example let $(w_i)_{i=1}^n$ be a normalized 1-unconditional 
1-subsymmetric finite basic sequence and let $\varep >0$. One can 
find a norm on $X$ such that every block basis contains a further block 
basis $(z_i)$ with $(z_{k_i})_{i=1}^n \buildrel {1+\varep}\over\sim 
(y_i)_{i=1}^n$ whenever $k_1<\cdots< k_n$. This is accomplished by 
taking (using the terminology of the proof of Theorem~3.1) 
$$\eqalign{
\Gamma & = \biggl\{ z^* = \sum_{i=1}^n b_i \sum_{j=(k_i-1)n+1}^{k_in} 
\!\!\!\! z_{m_j}^* : (z_{m_j}^*)_{j=1}^\infty\ \hbox{ is a block basis of }
\ (e_i^*)\cr
\myskip
&\qquad \hbox{ with }\ z_{m_1}^* \in 3Ba (X_1^*)\ ,\ 
z_{m_{j+1}}^* \in 3Ba (X_{m_j}^*)\ \hbox{ for }\ j\in \nat\ ,\cr
\myskip 
&\qquad k_1<k_2<\cdots < k_n\ \hbox{ and }\ 
\Big\| \sum_1^n b_i w_i^*\Big\| \le 1\biggr\}\ .\cr}$$ 

\proclaim Theorem 3.2. 
For $1<p<\infty$, $\ell_p$ is sequentially arbitrarily distortable. 

In order to prove Theorem 3.2 we will make use of the Banach space $S$   
introduced in [34].

The space $S$ 
has a 1-unconditional 1-subsymmetric normalized basis $(e_i)$ whose norm 
satisfies the following implicit equation 
$$\|x\| = \max \biggl\{ \|x\|_{c_0}\ ,\ 
\sup_{\scriptstyle \ell\ge2\atop\scriptstyle E_1<E_2<\cdots <E_\ell} 
{1\over \phi (\ell)} \sum_{i=1}^\ell \|E_i x\|\biggr\}$$ 
where $\phi (\ell) = \log_2 (1+\ell)$. 

The fact that $S$ is arbitrarily distortable [34] and complementably 
minimal [35] hinges heavily on two types of vectors which live in all 
block subspaces: $\ell_1^n+$ averages and averages of rapidly increasing 
$\ell_1^{n_i}+$ averages or RIS vectors. 
Precisely, following the terminology of [9], 
we call $x\in S$ an {\it $\ell_1^n+$ average with constant\/} 
$C$ if $\|x\|=1$ and $x= \sum_{i=1}^n x_i$ for some  
block basis $(x_i)_{i=1}^n$ 
of $(e_i)$ where $\|x_i\|\le Cn^{-1}$ for all $i$. 

Let $M_\phi (x) = \phi^{-1} (36x^2)$ for $x\in \real$. A block basis 
$(x_i)_{i=1}^N$ is an {\it RIS of length $N$ with constant\/} $C\equiv 1+
\varep <2$ if each $x_k$ is an $\ell_1^{n_k}+$ average with constant $C$, 
$n_1\ge 2CM_\phi (N/\varep)/2\varep \ln 2$ and $\frac{\varep}2 \phi (n_k)^
{1/2} \ge |\supp (x_{k-1})|$ for $k=2,\ldots,N$. The vector 
$x= \sum_{i=1}^N x_i/\|\sum_{i=1}^N x_i\|$ is called an 
{\it RIS vector of length $N$ and constant\/} $C$ and we say that the RIS 
sequence $(x_i)_{i=1}^N$ generates $x$.

\proclaim Lemma 3.3. {\rm [9]}\enspace 
Let $\varep_i\downarrow 0$. There exist integers $p_k\uparrow \infty$ and 
reals $\delta_k\downarrow 0$ so that if $A_k = \{x\in S: x$ is an RIS 
vector of length $p_k$ with constant $1+\delta_k\}$ 
and $A_k^* = \{x^* \in S^* :x^* = {1\over\phi (p_k)} 
\sum_1^{p_k} x_i^*$ where $(x_i^*)_1^{p_k}$ is a block sequence in 
$Ba(S^*)\}$ then 
\smallskip
\iitem{\rm a)} $|x_k^*(x_\ell)| <\varep_{\min (k,\ell)}$ if $k\ne\ell$, 
$x_k^* \in A_k^*$ and $x_\ell \in A_\ell$. 
\iitem{\rm b)} For all $k\in \nat$ and $x\in A_k$ there exists $x^*\in A_k^*$ 
with $x^*(x) > 1-\varep_k$.  In fact if $x$ is generated by 
$(x_i)_{i=1}^{p_k}$, $x^*$ may be taken to be ${1\over \phi (p_k)} \sum_{i=1}^
{p_k} x_i^*$ where $(x_i^*)_1^{p_k}$ is a normalized block basis of 
$(e_i^*)$ which is biorthogonal to $(x_i)_1^{p_k}$.
\smallskip
\noindent Moreover $A_k$ is asymptotic in $S$ for all $k\in\nat$. 

Using the sets $A_k$ and $A_k^*$ we can define the following 
subsets of $\ell_1$ 
$$B_k = \bigl\{ x_k^* \circ x_k/|x_k^*| (|x_k|) : x_k^*\in A_k^*\ ,\ 
x_k\in A_k \hbox{ and } |x_k^*|(|x_k|) = \|x_k^*\circ x_k\|_1 \ge 1-\varep_k
\bigr\}\ .$$ 

A set of sequences $B$ is 
{\it unconditional\/} if $x= (x_i)\in B$ implies that $(\pm x_i)\in B$ 
for all choices of signs and $B$ is {\it spreading\/} if $x= (x_i)\in B$ 
implies $\sum_i x_i e_{n_i} \in B$ for all increasing sequences $(n_i)$. 
Note that $A_k^*\subseteq Ba(S^*)$ and the sets $A_k$ and $A_k^*$ are 
unconditional and spreading. Thus the sets $B_k \subseteq S(\ell_1)$ are 
also spreading and unconditional. 

\proclaim Theorem 3.4. 
The sets $B_k\subseteq S(\ell_1)$, $k\in \nat$, are unconditional, 
spreading and asymptotic. 

We postpone the proof of Theorem 3.4.


\demo Proof of Theorem 3.2.
We first give the argument for $p=2$. Let $C_k = \{ v\in S(\ell_2) : 
|v|^2 \in B_k\}$. $C_k$ is just the  image of $B_k$ in $S(\ell_2)$ 
under the Mazur map. Since the Mazur map preserves block subspaces 
and is a uniform homeomorphism, $C_k$ is asymptotic in $\ell_2$ 
for all $k$. Moreover the $C_k$'s are nearly biorthogonal. Indeed 
if $v_k \in C_k$, $v_\ell\in C_\ell$ with $k\ne \ell$ let $|v_k|^2 = 
(x_k^*\circ x_k)/|x_k^*|(|x_k|)$ and $|v_\ell|^2 = 
(x_\ell^* \circ x_\ell)/|x_\ell^*|(|x_\ell|)$ be as in 
the definition of $B_k$ and $B_\ell$. Then 
letting $\lambda = (1-\varep_1)^{-1}$ 
$$\eqalign{
\langle |v_k|,|v_\ell|\rangle 
& \le \lambda  \sum_j \big| x_k^* (j)x_k(j)x_\ell^* (j)x_\ell(j)\big|^{1/2}\cr 
\myskip 
&\le \lambda\biggl( \sum_j \big| x_k^* (j)x_\ell(j)\big| \biggr)^{1/2} 
\biggl( \sum_j \big| x_\ell^* (j) x_k(j)\big| \biggr)^{1/2}\ 
\hbox{ (by Cauchy-Schwarz)}\cr 
\myskip 
& =\lambda \langle |x_k^*|, |x_\ell|\rangle^{1/2} 
\langle |x_\ell^*|,|x_k|\rangle^{1/2}
\le \lambda \varep_{\min (k,\ell)}\ \hbox{ (by Lemma 3.3).}\cr}$$ 
Define $\|x\|_k = \sup \{|\langle x,v\rangle| : v\in C_k \cup \varep_k 
Ba (\ell_2)\}$. 

If $p\ne 2$ we use a similar argument. 
Let $C_k = \{v\in S(\ell_p) : |v|^p\in B_k\}$ 
and $D_k = \{v\in S(\ell_q) : |v|^q\in B_k\}$ 
where ${1\over p}+{1\over q} =1$. Define $\|\cdot\|_k$ on $\ell_p$  by 
$$\|x\|_k = \sup \bigl\{ |\langle x,v\rangle| : v\in D_k \cup 
\varep_k Ba (\ell_q)\bigr\}\ .$$ 
Again, via the Mazur map, $C_k$ is asymptotic in $\ell_p$. 

Let $v_k\in C_k$ and $v_\ell\in D_\ell$ with $k\ne \ell$. 
Let $|v_k|^p = (x_k^*\circ x_k)/|x_k^*|(|x_k|)$ 
and $|v_\ell|^q = (x_\ell^*\circ x_\ell)/|x_\ell^*|(|x_\ell|)$ 
be as in the definition of $B_k$ and $B_\ell$. Assume $p>2$. Then 
$$\eqalign{|\langle |v_k|,|v_\ell|\rangle | 
& \le\lambda \sum_j \big| x_k^* (j) x_k(j)\big|^{1/p} 
\big| x_\ell^* (j) x_\ell (j)\big|^{1/q} \cr 
\myskip 
& = \lambda\sum_j \big| x_k^*(j) x_k(j) x_\ell^* (j) x_\ell (j)\big|^{1/p} 
\big| x_\ell^* (j) x_\ell(j)\bigr|^{{1\over q}-{1\over p}}\ .\cr}$$ 
Using H\"older's inequality with exponents $p\over 2$ and $p\over p-2$ 
and the fact that ${1\over q}-{1\over p} = {p-2\over p}$ we obtain  
that the last expression is 
$$ \eqalign{
&\le \lambda\biggl( \sum_j \big|x_k^*(j)x_k(j)x_\ell^*(j)x_\ell(j)\big|^{1/2}
\biggr)^{2/p} 
\biggl(\sum_j \big| x_\ell^*(j)x_\ell(j)\big| \biggr)^{p-2\over p} \cr 
\myskip 
&\le \lambda \varep_{\min(k,\ell)}^{2/p}\cr}$$ 
from the first part of the proof. The same estimates 
prevail if $p<2$.\qed 

\demo Remark. 
The proof yields that for $1<p<\infty$, ${1\over p} + {1\over q} =1$ there 
exist sequences $C_k \subseteq S(\ell_p)$ and $D_k \subseteq S(\ell_q)$ 
of nearly biorthogonal asymptotic unconditional spreading sets. 

It remains only to prove Theorem 3.4 which entails only showing that each 
$B_k$ is asymptotic. This will follow from the following 

\proclaim Lemma 3.5. 
Let $Y$ be a block  subspace of $\ell_1$ and let $\varep >0$, 
$m\in\nat$. There exists a vector $u\in S$ which is an $\ell_1^m+$ average 
with constant $1+\varep$ and $u^* \in Ba(S^*)$ with $d(u^*\circ u,S(Y)) 
<\varep$. 

Indeed assume that the lemma is proved and let $k\in \nat$ and $\varep>0$. 
>From the lemma we can find finite block sequences $(u_i)_{i=1}^{p_k} 
\subseteq S(S)$ and $(u_i^*)_{i=1}^{p_k}\subseteq Ba (S^*)$ along with a 
normalized block sequence $(y_i)_{i=1}^{p_k} \subseteq S(Y)$ such that 
\smallskip
\iitem{1)} $u = (\sum_{i=1}^{p_k} u_i)/\|\sum_{i=1}^{p_k} u_i\|$ 
is an RIS vector of length $p_k$ and constant $(1+\delta_k)$ 
generated by the RIS $(u_i)_{i=1}^{p_k}$. 
\iitem{2)} $\|u_i^*\circ u_i - y_i\|_1 < \varep$ for $i\le p_k$. 
\iitem{3)} $u_i^*\circ u_j =0$ if $i\ne j$.
\smallskip
\noindent Let $u^* = {1/\phi(p_k)}\sum_1^{p_k} u_i^*$. Then $u^*\in A_k^*$ 
and $(u^*\circ u/\|u^*\circ u\|_1) \in B_k$ by Lemma 3.3~b). 
Now 
$$\|u^* \circ u\|_1 = {p_k\over \phi (p_k) \|\sum_1^{p_k} u_i\|}$$ 
so 
$${u^*\circ u\over \|u^*\circ u \|_1} = {1\over p_k} \sum_{i=1}^{p_k} 
u_i^*\circ u_i\ .$$ 
Thus 
$$\Big\| {u^*\circ u\over \|u^*\circ u\|_1} - {1\over p_k} 
\sum_{i=1}^{p_k} y_i\Big\|_1 \le {1\over p_k} \sum_{i=1}^{p_k} 
\|u_i^*\circ u_i-y_i\|_1 
<\varep\quad\hbox{by 2).}$$

This proves that $B_k$ is asymptotic in $\ell_1$. 

In order to prove Lemma 3.5 we first need a sublemma. We denote the maps 
$E_{S^*}(h)$ and $F_{S^*}(h)$ by  $E_*(h)$ and $F_*(h)$, respectively.  

\proclaim Sublemma 3.6. 
Let $m,K$ be integers and let $0<\tau <1$ be such that $\log \phi (m^K) < 
\tau K$. Let $(h_i)_{i=1}^{m^K}$ be a normalized block sequence in 
$\ell_1^+$. Then there exist in $\ell_1^+$ a normalized block basis 
$(b_i)_{i=1}^m$ of $(h_i)_1^{m^K}$ such that 
$$\sum_{j=1}^m E_* (b_j) - E_* \biggl( \sum_{j=1}^m b_j\biggr) \le \tau m \ .
\eqno(3.1)$$ 

\proof 
For each $i\le m^K$, let $v_i = F_*(h_i)$. 
Now ${1\over \phi (m^K)} \sum_1^{m^K} v_i \in Ba(S^*)$ and so 
$$\eqalignno{
E_* \biggl(\sum_1^{m^K} h_i\biggr) 
& \ge E\biggl( \sum_1^{m^K} h_i,{1\over \phi (m^K)} \sum_1^{m^K} v_i\biggr) 
&(3.2)\cr 
\myskip 
&= \sum_1^{m^K} E(h_i,v_i) - m^K \log \phi (m^K)\cr 
\myskip 
& = \sum_1^{m^K} E_* (h_i) - m^K \log \phi (m^K)\ .\cr}$$ 
Let $\sum_{i=1}^{m^K} h_i = \sum_{j=1}^m d_j^1$     where 
$(d_j^1)_{j=1}^m$ is a block basis of $(h_i)$, each $d_j^1$ consisting of 
the sum of $m^{K-1}$ of the $h_i$'s. Break each $d_j^1$ into $m$ 
successive pieces, each containing $m^{K-2}$ of the $h_i$'s to obtain 
$d_j^1 = \sum_{\ell=1}^m d_{j,\ell}^2$ and continue to define 
$d_{\alpha,j}^\ell$ for $\ell \le k$ and $\alpha\in \{1,\ldots,m\}^{\ell-1}$ 
in this fashion. 
Consider the telescoping sum 
$$\eqalign{ 
\sum_{i=1}^{m^K} E_* (h_i) - E_*\biggl( \sum_{i=1}^{m^K} h_i\biggr) 
& = \sum_{j=1}^m E_* (d_j^1) - E_* \biggl(\sum_{j=1}^m d_j^1\biggr) \cr 
\myskip 
&\qquad + \sum_{j=1}^m \left[ \sum_{\ell=1}^m E_* (d_{j,\ell}^2) 
- E_* \biggl( \sum_{\ell=1}^m d_{j,\ell}^2\biggr) \right] 
+ \cdots \ .\cr}$$ 

For $1\le s\le K$, the $s^{th}$ level of this decomposition is the sum of 
$m^{s-1}$ nonnegative terms of the form  (for $\alpha\in\{1,\ldots,m\}^{s-1}$) 
$$\sum_{\ell=1}^m E_* (d_{\alpha,\ell}^s) 
- E_*\biggl( \sum_{\ell=1}^m d_{\alpha,\ell}^s\biggr)\ . 
\eqno(3.3)$$ 
If each of these terms is greater than $\tau m^{K-s+1}$ then the sum of all  
terms on the $s^{th}$ level is greater than $\tau m^K$ and so the sum 
over all $K$ levels yields 
$$\sum_{i=1}^{m^K} E_* (h_i) - E_* \biggl( \sum_1^{m^K} h_i\biggr) 
> K\tau m^K$$ 
which contradicts (3.2). 

Thus the number (3.3) does not exceed the value 
$\tau m^{K-s+1}$ for some $s$ and multi-index $\alpha$. 
Let $b_\ell = d_{\alpha,\ell}^s /\|d_{\alpha,\ell}^s\|$. 
Using  $E_* (ah) = aE_* (h)$ for $a>0$ and $\|d_{\alpha,\ell}^s\| 
= m^{K-s}$ we obtain 
$$\sum_{\ell=1}^m E_* (b_\ell) - E_* \biggl( \sum_{\ell=1}^m b_\ell\biggr) 
\le {\tau m^{K-s+1}\over m^{K-s}} = \tau m\ .\eqno\blackbox$$ 

\demo Proof of Lemma 3.5. 
Let $\varep >0$, $m\in \nat$ and let $Y$ be a block subspace of $\ell_1$ 
with block basis $(h_i)$. 
By unconditionality in $S$ it suffices to consider only the case where 
$(h_i) \subseteq S(\ell_1)^+$. Let $0<\tau < {\psi (\varep)
\over m}$ (see Definition~2.2) and choose $K\in \nat$ such that 
$\tau K >\log (\phi (m^K))$. By sublemma~3.6 choose a block basis 
$(b_i)_1^m$ of $(h_i)_{i=1}^{m^K}$, $(b_i)_1^m\subseteq S(\ell_1^+)$ with 
$$\sum_1^m E_* (b_j) - E_* \biggl( \sum_{j=1}^m b_j\biggr) < \tau m\ . 
\eqno(3.4)$$ 
Choose $x^*= F_*(\sum_{j=1}^m b_j)$ 
and write $x^* =\sum_{j=1}^m x_j^*$ 
with $\supp x_j^* = \supp b_j$. For $j\le m$ let $w_j^* = F_* (b_j)$.  
As we noted in Section~2, for each $j$ 
there exists $w_j\in S(S)^+$ with $b_j = w_j^*\circ w_j$ and  
$\supp w_j = \supp b_j$. By (3.4) we have 
$$\sum_{j=1}^m E(b_j,w_j^*) - E\biggl( \sum_{j=1}^m b_j,x^*\biggr) 
= \sum_{j=1}^m \bigl[ E(b_j,w_j^*) - E(b_j,x_j^*)\bigr] 
< \tau m < \psi (\varep)\ .$$ 
Since each term in the middle expression is nonnegative we obtain 
$$E(b_j,x_j^*) > E(b_j,w_j^*) - \psi (\varep)\ \hbox{ for }\ 
j\le m\ .$$ 
By Proposition 2.3(A) there exists sets $H_j\subseteq \supp b_j$ such that 
$\|H_jb_j\|_1 > 1-\varep$ and 
\vskip1pt
\noindent $(1-\varep) H_j w_j^* \le H_jx_j^* 
\le (1+\varep) H_jw_j^*$ pointwise for all $1\le j\le m$. 

$H_jb_j = H_jw_j^* \circ w_j$ and 
$\|H_jx_j^* - H_jw_j^*\|\le \varep$ so $\|H_jb_j - H_jx_j^*\circ w_j\|_1 
\le \varep$. Thus 
$$\|b_j - H_jx_j^* \circ w_j\|_1 \le 2\varep\ ,\ \hbox{ for }\ 
1\le j\le m\ .\eqno(3.5)$$ 
>From this we first note that $H_jx_j^*(w_j) \ge 1-2\varep$ and so for 
$a_i$'s nonnegative, 
$$\eqalign{ 
\Big\| \sum_1^m a_jw_j\Big\| & \ge x^*\biggl( \sum_1^m a_j w_j\biggr) \cr 
\myskip 
&\ge \sum_{j=1}^m a_j H_j x_j^*(w_j) \ge \biggl( \sum_{j=1}^m a_j\biggr) 
(1-2\varep)\ .\cr}$$ 
By unconditionality $(w_j)_{j=1}^m$ is an $\ell_1^m$ sequence with  
constant $(1-2\varep)^{-1}$. 

Secondly, set 
$$\bar w = {1\over m} \sum_{j=1}^m w_j\quad \hbox{and}\quad 
w= {1\over \big\|\sum\limits_{j=1}^m w_j\big\|} \ 
\sum_{j=1}^m w_j\ .$$ 
$w$ is an $\ell_1^m$ average with constant $(1-2\varep)^{-1}$. Furthermore 
$$\eqalign{\Big\| {1\over m} \sum_{j=1}^m b_j - \biggl( \bigcup_{j=1}^m 
H_j\biggr) x^*\circ w\Big\|_1 
&\le  \Big\| {1\over m} \sum_{j=1}^m b_j - 
\biggl( \bigcup_{j=1}^n H_j\biggr) x^*\circ \bar w\Big\|_1 
+ \|w-\bar w\| \cr 
\myskip 
&\le {1\over m} \sum_{j=1}^m \| b_j - H_j x^*\circ w_j\|_1 
+ \|w-\bar w\|\ .\cr}$$ 
The first term is $< 2\varep$ by (3.5). 
Since $\|\sum_{j=1}^m w_j\| \ge m(1-2\varep)$, 
$\|w-\bar w\| \le 2\varep /1-2\varep$. 
Thus 
$$ d\left( \biggl( \bigcup_{j=1}^m H_j\biggr) x^* \circ w\ ,\ S(Y)\right) 
< 2\varep + {2\varep\over 1-2\varep}$$ 
which proves Lemma 3.5.\qed 

\demo Remark 3.7. 
B.~Maurey [22] has recently extended the results above. He has proven that 
if $X$ has an unconditional basis and is superreflexive,   
then $X$ contains an arbitrarily distortable subspace. 
He has also used a modification of this argument to give a simpler proof 
of Milman's result that every Banach space contains $c_0$ or $\ell_p$ 
($1\le p<\infty$) or a distortable subspace. 
\vskip1pt
N.~Tomczak-Jaegermann and V.~Milman [29] have proven that if $X$ 
has bounded distortion, then $X$ contains an ``asymptotic $\ell_p$ 
or $c_0$.'' $X$ has bounded distortion if for some $\lambda <\infty$, no 
subspace of $X$ is $\lambda$-distortable. A space with a basis $(e_i)$ 
is an asymptotic $\ell_p$ if for some $C<\infty$ for all $n$ whenever 
$$e_n < x_1 <\cdots < x_n\ ,\qquad \|x_i \| =1\quad (i=1,\ldots,n)\ ,$$ 
then $(x_i)_1^n$ is $C$-equivalent to the unit vector basis of $\ell_p^n$. 
\bigskip 

\centerline{\cbf References}
\medskip

{\ninepoint

\item  {[1]} Y. Benyamini, 
{The uniform classification of Banach spaces}, 
{\it Longhorn Notes 1984-85},
The University of Texas at Austin, 1985. 

\item  {[2]} F. Chaatit, 
{Uniform homeomorphisms between unit spheres of Banach lattices}, 
in preparation. 

\item  {[3]} J. Diestel, 
{\it Geometry of Banach Spaces -- Selected Topics},
LNM 485, Springer-Verlag, New York, 1975. 

\item  {[4]} P. Enflo, 
{On a problem of Smirnov}, 
{\it Ark. Mat.} {\bf8} (1969), 107--109. 

\item  {[5]} T. Figiel and W.B. Johnson, 
{A uniformly convex Banach space which contains no $\ell_p$}, 
{\it Comp. Math.} {\bf29} (1974), 179--190.

\item  {[6]} T.A. Gillespie, 
{Factorization in Banach function spaces}, 
{\it Indagationes Math.} {\bf43} (1981), 287--300. 

\item  {[7]} T. Gowers, 
{Lipschitz functions on classical spaces},  
{\it European J. Comb.}, to appear. 

\item  {[8]} T. Gowers, 
Ph.D.\ Thesis, Cambridge University, 1990. 

\item  {[9]} T. Gowers and B. Maurey, 
{The unconditional basic sequence problem}, 

\item  {[10]} R. Haydon, E. Odell, H. Rosenthal and Th. Schlumprecht, 
{On distorted norms in Banach spaces and the existence of $\ell_p$ types}, 
preprint. 

\item  {[11]} R.C. James, 
{Uniformly nonsquare Banach spaces}, 
{\it Ann. Math.} {\bf80} (1964), 542--550.

\item  {[12]} R.C. James, 
{A nonreflexive Banach space that is uniformly nonoctohedral}, 
{\it Israel J. Math.} {\bf18} (1974), 145--155. 

\item  {[13]} R.E. Jamison and W.H. Ruckle, 
{Factoring absolutely convergent series}, 
{\it Math. Ann.} {\bf224} (1976), 143--148. 

\item  {[14]} N.J. Kalton, 
{Differentials of complex interpolation processes for Kothe function 
spaces}, preprint. 

\item  {[15]} N.J. Kalton, 
{Uniform homeomorphisms and complex interpolation}, 
in preparation. 

\item  {[16]} J.L. Krivine, 
{Sous espaces de dimension finie des espaces de Banach r\'eticul\'es}, 
{\it Ann. of Math.} {\bf104} (1976), 1--29. 

\item  {[17]} J. Lindenstrauss and A. Pe{\l}czy\'nski, 
{Contributions to the theory of classical Banach spaces}, 
{\it J.~Funct. Anal.} {\bf8} (1971), 225--249. 

\item  {[18]} J. Lindenstrauss and L. Tzafriri, 
{\it Classical Banach Spaces, I},
Springer-Verlag, New York, 1977. 

\item  {[19]} J. Lindenstrauss and L. Tzafriri, 
{\it Classical Banach Spaces, II},
Springer-Verlag, New York, 1979. 

\item  {[20]} G.Ya. Lozanovskii, 
{On some Banach lattices}, 
{\it Sibir. Math. J.} {\bf10} (1969), 584--99. 

\item  {[21]} G.Ya. Lozanovskii, 
{On some Banach lattices III}, 
{\it Sibir. Math. J.} {\bf13} (1972), 1304--13. 

\item  {[22]} B. Maurey, 
{A remark about distortion}, 
in preparation. 

\item  {[23]} B. Maurey and G. Pisier, 
{S\'eries de variables al\'eatoires vectorielles ind\'ependantes et 
propri\'et\'es g\'eom\'etriques des espaces de Banach}, 
{\it Studia Math.} {\bf58} (1976), 45--90. 

\item  {[24]} B. Maurey and H. Rosenthal, 
{Normalized weakly null sequences with no unconditional subsequences}, 
{\it Studia Math.} {\bf61} (1977), 77--98. 

\item  {[25]} S. Mazur, 
{Une remarque sur l'hom\'eomorphisme des champs fonctionnels}, 
{\it Studia Math.} {\bf1} (1930), 83--85. 

\item  {[26]} V. Milman, 
{Geometric theory of Banach spaces II, geometry of the unit sphere}, 
{\it Russian Math Survey} {\bf26} (1971), 79--163 (trans. from Russian). 

\item  {[27]} V. Milman, 
{The spectrum of bounded continuous functions defined on the unit 
sphere of a $B$-space}, 
{\it Functsional Analiz. i ego Prilozhen} {\bf3} (1967), 67--79. 

\item  {[28]} V. Milman and G. Schechtman, 
{\it Asymptotic Theory of Finite Dimensional Normed Spaces},
LNM 1200, Springer-Verlag, New York, 1986. 

\item  {[29]} V. Milman and N. Tomczak-Jaegermann, 
{Asymptotic $\ell_p$ spaces and bounded distortions}, 
preprint. 

\item  {[30]} G. Pisier, 
{The volume of convex bodies and Banach space geometry}, 
{\it Cambridge Tracts in Math.} {\bf94} (1989). 

\item  {[31]} Y. Raynaud, 
{Espaces de Banach superstables, distances stables et 
hom\'eomorphismes uniformes}, 
{\it Israel J. Math.} {\bf44} (1983), 33--52. 

\item  {[32]} M. Ribe, 
{Existence of separable uniformly homeomorphic non isomorphic 
Banach spaces}, 
{\it Israel J. Math.} {\bf48} (1984), 139--147. 

\item  {[33]} H. Rosenthal, 
{A characterization of Banach spaces containing $\ell_1$}, 
{\it Proc. Nat. Acad. Sci.} (U.S.A.) {\bf71} (1974), 2411--2413. 

\item  {[34]} Th. Schlumprecht, 
{An arbitrarily distortable Banach space},

\item  {[35]} Th. Schlumprecht, 
{A complementably minimal Banach space not containing $c_0$ or $\ell_p$}, 

\item  {[36]} S. Szarek and N. Tomczak-Jaegermann, 
{On nearly euclidean decomposition for some classes of Banach spaces}, 
{\it Comp. Math.} {\bf40} (1980), 367--385. 

\item  {[37]} Nicole Tomczak-Jaegermann, 
{Banach-Mazur distances and finite-dimensional operator ideals}, 
{\it Pitman Monographs} {\bf38} (1989). 

\item  {[38]}  B.S. Tsirelson,
{Not every Banach space contains $\ell_p$ or $c_0$}, 
{\it Funct. Anal. Appl.} {\bf 8} (1974), 138--141.

\bigskip
\baselineskip=12pt
{\vbox{\halign{\indent #\hfill\qquad &#\hfill\cr
E. Odell&Th. Schlumprecht \cr
Department of Mathematics&Department of Mathematics \cr
The University of Texas at Austin&Louisiana State University \cr
Austin, TX 78712-1082&Baton Rouge, LA 70803\cr}}}
\medskip 

{\obeylines
{\sl Current Address for Th. Schlumprecht:}
Department of Mathematics
Texas A\&M University 
College Station, TX 77843 
}
\medskip
\rightline{\sl September 8, 1992}

}
\end